# Symmetric Dyck paths and q-Narayana numbers

Johann Cigler

**Abstract**

We show that the q-Narayana numbers $N_{n,k}(q)$ for $q=-1$ can be interpreted as the number of symmetric Dyck paths of semi-length $n$ with $k$ valleys.

## 1. q-Narayana numbers as Hoggatt binomials

In addition to several other combinatorial interpretations the Narayana numbers

(1) $$N_{n,k} = \frac{1}{k+1}\binom{n}{k}\binom{n-1}{k} = \frac{1}{n}\binom{n}{k}\binom{n}{k+1}$$

count Dyck paths of semi-length $n$ with $k$ valleys (cf. e.g. [5]). Algebraically they can be represented as so-called Hoggatt binomials (cf. [1],[2]). To this end let

(2) $$\langle n \rangle := \binom{n+1}{2}, \quad \langle n \rangle! := \prod_{j=1}^{n}\langle j \rangle$$

and define

(3) $$\left\langle {n \atop k} \right\rangle = \frac{\langle n \rangle!}{\langle k \rangle!\langle n-k \rangle!} = \frac{\langle n \rangle}{\langle k \rangle}\left\langle {n-1 \atop k-1} \right\rangle$$

for $0 \le k \le n$ and $\left\langle {n \atop k} \right\rangle = 0$ else.

Then we get

(4) $$\left\langle {n \atop k} \right\rangle = \frac{1}{k+1}\binom{n}{k}\binom{n+1}{k} = N_{n+1,k}.$$

The first terms of $\left\langle {n \atop k} \right\rangle$ are

(5)

| n \ k | 0 | 1 | 2 | 3 | 4 | 5 | 6 | 7 |
|---|---|---|---|---|---|---|---|---|
| 0 | 1 | | | | | | | |
| 1 | 1 | 1 | | | | | | |
| 2 | 1 | 3 | 1 | | | | | |
| 3 | 1 | 6 | 6 | 1 | | | | |
| 4 | 1 | 10 | 20 | 10 | 1 | | | |
| 5 | 1 | 15 | 50 | 50 | 15 | 1 | | |
| 6 | 1 | 21 | 105 | 175 | 105 | 21 | 1 | |
| 7 | 1 | 28 | 196 | 490 | 490 | 196 | 28 | 1 |

For the $q-$Narayana numbers (cf. [3],[5])

(6) $$N_{n,k}(q) = \frac{1}{[k+1]}\begin{bmatrix}n\\k\end{bmatrix}\begin{bmatrix}n-1\\k\end{bmatrix}$$

we get an analogous representation by replacing the binomial coefficients with the $q-$binomial coefficients $\begin{bmatrix}n\\k\end{bmatrix} = \begin{bmatrix}n\\k\end{bmatrix}_q = \frac{[n]!}{[k]![n-k]!}$ where $[n] = [n]_q = \frac{1-q^n}{1-q}$ and $[n]! = \prod_{j=1}^{n}[j]$.

In this case we choose

(7) $$\langle n\rangle_q := \begin{bmatrix}n+1\\2\end{bmatrix}_q.$$

Then $\langle n\rangle_q! = \frac{[n]![n+1]!}{[2]^n}$ and $\frac{\langle n\rangle_q!}{\langle k\rangle_q!\langle n-k\rangle_q!} = \frac{[n]![n+1]!}{[k]![k+1]![n-k]![n+1-k]!}$

gives

(8) $$\left\langle\begin{matrix}n\\k\end{matrix}\right\rangle_q = N_{n+1,k}(q) = \frac{1}{[k+1]}\begin{bmatrix}n\\k\end{bmatrix}\begin{bmatrix}n+1\\k\end{bmatrix}.$$

For $q=-1$ we get

(9) $$\langle n\rangle_{-1} = \begin{bmatrix}n+1\\2\end{bmatrix}_{-1} = \lim_{q\to -1}\frac{[n][n+1]}{(1+q)} = \left\lfloor\frac{n+1}{2}\right\rfloor,$$

because for odd $n$ we get $[n]_{-1} = 1$ and for even $n = 2m$

$$\lim_{q\to -1}\frac{[n]}{(1+q)} = \lim_{q\to -1}\frac{[2m]}{(1+q)} = \lim_{q\to -1}\frac{1-q^{2m}}{1-q^2} = \frac{2m}{2} = m.$$

By (3) and (9) we get $\left\langle\begin{matrix}n\\k\end{matrix}\right\rangle_{-1} = \frac{\left\lfloor\frac{n+1}{2}\right\rfloor}{\left\lfloor\frac{k+1}{2}\right\rfloor}\left\langle\begin{matrix}n-1\\k-1\end{matrix}\right\rangle_{-1}$. Since $\frac{\left\lfloor\frac{n+1}{2}\right\rfloor}{\left\lfloor\frac{k+1}{2}\right\rfloor}\left(\begin{matrix}\left\lfloor\frac{n-1}{2}\right\rfloor\\\left\lfloor\frac{k-1}{2}\right\rfloor\end{matrix}\right) = \left(\begin{matrix}\left\lfloor\frac{n+1}{2}\right\rfloor\\\left\lfloor\frac{k+1}{2}\right\rfloor\end{matrix}\right)$ we get

by induction

(10) $$\left\langle\begin{matrix}n\\k\end{matrix}\right\rangle_{-1} = \left(\begin{matrix}\left\lfloor\frac{n}{2}\right\rfloor\\\left\lfloor\frac{k}{2}\right\rfloor\end{matrix}\right)\left(\begin{matrix}\left\lfloor\frac{n+1}{2}\right\rfloor\\\left\lfloor\frac{k+1}{2}\right\rfloor\end{matrix}\right).$$

The first terms of $\left\langle\begin{matrix}n\\k\end{matrix}\right\rangle_{-1}$ are

(11)

| n\k | 0 | 1 | 2 | 3 | 4 | 5 | 6 | 7 |
|---|---|---|---|---|---|---|---|---|
| 0 | 1 | | | | | | | |
| 1 | 1 | 1 | | | | | | |
| 2 | 1 | 1 | 1 | | | | | |
| 3 | 1 | 2 | 2 | 1 | | | | |
| 4 | 1 | 2 | 4 | 2 | 1 | | | |
| 5 | 1 | 3 | 6 | 6 | 3 | 1 | | |
| 6 | 1 | 3 | 9 | 9 | 9 | 3 | 1 | |
| 7 | 1 | 4 | 12 | 18 | 18 | 12 | 4 | 1 |

It turns out that the Hoggatt binomials $\left\langle {n \atop k} \right\rangle_{-1}$ satisfy simple recurrences. For $k \geq 1$ we get

(12) $$\left\langle {n \atop 2k} \right\rangle_{-1} = \left\langle {n-1 \atop 2k} \right\rangle_{-1} + \left\langle {n-1 \atop 2k-1} \right\rangle_{-1}$$

and for $k \geq 0$

(13) $$\left\langle {2n+1 \atop 2k+1} \right\rangle_{-1} = \left\langle {2n \atop 2k+1} \right\rangle_{-1} + \left\langle {2n \atop 2k} \right\rangle_{-1},$$

$$\left\langle {2n \atop 2k+1} \right\rangle_{-1} = \left\langle {2n-1 \atop 2k+1} \right\rangle_{-1} + \frac{k}{k+1}\left\langle {2n-1 \atop 2k} \right\rangle_{-1} = \left\langle {2n-1 \atop 2k+1} \right\rangle_{-1} + \left\langle {2n-1 \atop 2k} \right\rangle_{-1} - N_{n,k}.$$

Using (3) identity (12) can be written as $\frac{\langle n \rangle_{-1}}{\langle 2k \rangle_{-1}} \left\langle {n-1 \atop 2k-1} \right\rangle_{-1} = \frac{\langle n-2k \rangle_{-1}}{\langle 2k \rangle_{-1}} \left\langle {n-1 \atop 2k-1} \right\rangle_{-1} + \left\langle {n-1 \atop 2k-1} \right\rangle_{-1}$ and follows from $\langle n \rangle_{-1} = \langle n-2k \rangle_{-1} + \langle 2k \rangle_{-1}$ which reduces to $\left\lfloor \frac{n+1}{2} \right\rfloor = \left\lfloor \frac{n+1-2k}{2} \right\rfloor + \left\lfloor \frac{2k}{2} \right\rfloor$.

Identity (13) follows in the same way. For the last identity we must show that $\frac{1}{k+1}\left\langle {2n-1 \atop 2k} \right\rangle_{-1} = N_{n,k}$. This follows from (10) which gives

$$\frac{1}{k+1}\left\langle {2n-1 \atop 2k} \right\rangle_{-1} = \frac{1}{k+1}\left(\left\lfloor \frac{2n-1}{2} \right\rfloor \atop \left\lfloor \frac{2k}{2} \right\rfloor\right)\left(\left\lfloor \frac{2n}{2} \right\rfloor \atop \left\lfloor \frac{2k+1}{2} \right\rfloor\right) = \frac{1}{k+1}\binom{n-1}{k}\binom{n}{k} = N_{n,k}.$$

Therefore the table $\left(\left\langle {n \atop k} \right\rangle_{-1}\right)_{n,k \geq 0}$ is uniquely determined by the recurrences (12) and (13) and the initial values $\left\langle {0 \atop k} \right\rangle = [k=0]$ and $\left\langle {n \atop 0} \right\rangle = 1$.

Let us also note that

(14) $$s(n) = \sum_{k=0}^{n} \left\langle {n \atop k} \right\rangle_{-1} = \binom{n+1}{\left\lfloor \frac{n+1}{2} \right\rfloor}.$$

Let $C_n = \frac{1}{n+1}\binom{2n}{n}$ denote the Catalan numbers. It is easily verified that $a(n) = \binom{n+1}{\left\lfloor \frac{n+1}{2} \right\rfloor}$ satisfies $a(2n+1) = 2a(2n)$ and $a(2n) = 2a(2n-1) - C_n$. We must show that $s(n)$ satisfies the same recurrences. This follows immediately from (12) and (13) and from $\sum_{k=0}^{n} N_{n,k} = C_n$.

## 2. Symmetric Dyck paths

The first terms of table (11) lead to OEIS [4], A088855, which suggests that $\left\langle {n-1 \atop k} \right\rangle_{-1}$ can be interpreted as the number of symmetric Dyck paths of semi-length $n$ with $k$ valleys. We now give a proof of this fact.

**Theorem 1**

Let $u(n,k)$ be the number of symmetric Dyck paths of semi-length $n$ with $k$ valleys and $v(n,k) = \left\langle {n-1 \atop k} \right\rangle_{-1} = N_{n,k}(-1)$.

Then

(15) $$u(n,k) = v(n,k) = \binom{\left\lfloor \frac{n-1}{2} \right\rfloor}{\left\lfloor \frac{k}{2} \right\rfloor}\binom{\left\lfloor \frac{n}{2} \right\rfloor}{\left\lfloor \frac{k+1}{2} \right\rfloor}.$$

**Proof**

A Dyck path $p$ of semi-length $n$ is a lattice path from $(0,0)$ to $(2n,0)$, consisting of $n$ up steps $U = (1,1)$ and $n$ down steps $D = (1,-1)$, such that the path never goes below the $x$-axis. We write more precisely $U = U_i$ if $i$ is the $x$-coordinate of the initial point of $U$ and $D = D_i$ if $i$ is the $x$-coordinate of the end point of $D$. A point $(i, p_i)$ defined by consecutive steps $D_i U_i$ is called a valley.

If we describe $p$ by its heights $p_i$ at position $i$ then a symmetric Dyck path with semi-length $n$ is a Dyck path $p = (0, p_1, p_2, \cdots, p_{2n-1}, 0)$ satisfying $p_i = p_{2n-i}$ for $0 \leq i \leq n$. Such a path $p$ is uniquely determined by the initial part $(0, p_1, \cdots, p_n)$ from $(0,0)$ to some $(n, j)$.

The following table shows the smallest symmetric Dyck paths classified according to the number $k$ of their valleys.

| $n \backslash k$ | 0 | 1 | 2 | 3 |
|---|---|---|---|---|
| 1 | $UD$ | | | |
| 2 | $U^2D^2$ | $UDUD$ | | |
| 3 | $U^3D^3$ | $U^2DUD^2$ | $UDUDUD$ | |
| 4 | $U^4D^4$ | $\{U^3DUD^3, U^2D^2U^2D^2\}$ | $\{U^2DUDUD^2, UDU^2D^2UD\}$ | $UDUDUDUD$ |

(16)

Let $S_n$ denote the set of non-negative paths of length $n$ starting at $(0,0)$ and let $G(n,j,k)$ be the number of paths $p \in S_n$ with $k$ valleys $DU$ and with end point $(n,j)$, where the last step is an up step $U$ and $g(n,j,k)$ be those paths where the last step is a down step $D$.

Since a down step at the last position $n$ defines a valley in the symmetric Dyck path with semi-length $n$ we get

(17)
$$u(n, 2k) = \sum_{j=0}^{n} G(n, j, k),$$
$$u(n, 2k+1) = \sum_{j=0}^{n} g(n, j, k).$$

It should be noted that only points $(n,j)$ with $j \equiv n \pmod 2$ can appear as endpoints. For other $(n,j)$ we have $G(n,j,k) = g(n,j,k) = 0$. For $j < 0$ we set $G(n,j,k) = g(n,j,k) = 0$. Note also that $G(n,0,k) = 0$.

For $k \geq 1$ we get

(18)
$$G(n, j, k) = G(n-1, j-1, k) + g(n-1, j-1, k-1).$$

Note that a path $pU : (0,0) \to (n, j)$ with $k$ valleys is of the form $pU = p_1UU + p_2DU$ where $p_1U$ has $k$ valleys and $p_2D$ has $k-1$ valleys because its last down step produces a valley of $pU$.

For paths ending with a down step we get

(19)
$$g(n, j, k) = g(n-1, j+1, k) + G(n-1, j+1, k).$$

Let us further note that $G(n,j,k) = 0$ for $n < 2k+1$ and that $G(2k+1, j, k) = [j = 1]$ because the smallest path with $k$ valleys and ending with $U$ is $(UD)^k U$ from $(0,0)$ to $(2k+1, 1)$.

Similarly $g(n,j,k) = 0$ for $n < 2k+2$ and $g(2k+2, j, k) = [j = 0]$ because the smallest path with $k$ valleys and ending with a down step is $(UD)^{k+1}$ from $(0,0)$ to $(2k+2, 0)$.

With these initial values $G(n,j,k)$ and $g(n,j,k)$ are uniquely determined.

By (17), (18) and (19) we get

(20)
$$u(n, 2k) = u(n-1, 2k) + u(n-1, 2k-1).$$

(21) $$u(2n, 2k+1) = u(2n-1, 2k+1) + u(2n-1, 2k)$$

and

(22) $$u(2n+1, 2k+1) = u(2n, 2k+1) + u(2n, 2k) - g(2n, 0, k).$$

The last term occurs because

$$\sum_{j=0}^{n} g(2n, j+1, k) = \sum_{j=0}^{n} g(2n, j, k) - g(2n, 0, k) = u(2n, 2k) - g(2n, 0, k).$$

Note that $g(2n, 0, k)$ is the number of Dyck paths of semi-length $n$ with $k$ valleys which is known to be $N_{n,k}$.

Thus $u(n,k)$ satisfies the same recurrences as $v(n,k)$. To prove that $u(n,k) = v(n,k)$ we need only verify that $u(n,0) = v(n,0) = 1$ and $u(n,1) = v(n,1) = \left\lfloor \frac{n+1}{2} \right\rfloor$.

For $k = 0$ the only non-negative path of length $n$ without valleys which ends with an up step is $U^n$.

The only paths of length $n$ without valleys which end with a down step are the paths
$U^{n-j} D^j : (0,0) \to (n, n-2j), \ 1 \leq j \leq \left\lfloor \frac{n}{2} \right\rfloor.$

Therefore

(23) $$u(n,1) = \left\lfloor \frac{n}{2} \right\rfloor = v(n,1) = \left\langle \begin{matrix} n-1 \\ 1 \end{matrix} \right\rangle.$$

## 3. Final remarks

If Dyck paths are counted according to major index then $q^{k^2+k} N_{n,k}(q)$ counts the paths with $k$ valleys (cf. [3],[5]). More precisely

(24) $$\sum_{k=0}^{n} q^{k^2+k} N_{n,k}(q) = C_n(q) = \frac{1}{[n+1]} \begin{bmatrix} 2n \\ n \end{bmatrix}.$$

For $q = -1$ this identity reduces to

(25) $$\sum_{k=0}^{n} N_{n,k}(-1) = \left( \begin{matrix} n \\ \left\lfloor \frac{n}{2} \right\rfloor \end{matrix} \right).$$

As shown above $\binom{n}{\left\lfloor \frac{n}{2} \right\rfloor}$ can be interpreted as the number of symmetric Dyck paths of semi-length $n$ and $N_{n,k}(-1)$ as the number of those paths with $k$ valleys.

It would be interesting to know whether this interpretation of (25) could also be deduced from the above interpretation of (24) or whether it just happens to look similar by pure numerical coincidence.